\documentclass[a4paper, 11pt, reqno]{amsart}
\usepackage{amsmath}
\usepackage{amsfonts}
\usepackage{amssymb}
\usepackage{amsthm}
\usepackage{mathtools}
\usepackage{hyperref}
\usepackage{indentfirst}
\usepackage[english]{babel}
\usepackage[T1]{fontenc}
\usepackage{graphicx}
\usepackage{tensor} 
\usepackage[colorinlistoftodos]{todonotes}

\usepackage{subfiles}

\newcommand{\si}{\Sigma}

\newcommand{\tr}{\operatorname{tr}}
\newcommand{\diver}{\operatorname{div}}

\DeclareMathOperator{\Ric}{Ric}

\newcommand{\Ein}{\operatorname{G}}

\newcommand{\s}{\operatorname{Sc}}

\newtheorem{theo}{Theorem}[]

\newtheorem{lemma}[theo]{Lemma}
\newtheorem{definition}[theo]{Definition}
\newtheorem{coro}[theo]{Corollary}

\newtheorem{remark}{Remark}
\newcommand*\Item[1][]{%
  \ifx\relax#1\relax  \item \else \item[#1] \fi
  \abovedisplayskip=0pt\abovedisplayshortskip=0pt~\vspace*{-\baselineskip}
}

\begin{document}

\title[THE KERR-NEWMAN-DE SITTER SPACETIME] 
{The index of the cosmological horizon and the area-charge-inequality} 


\author{Neilha Pinheiro}
\address{
Universidade Federal do Amazonas\\
Instituto de Ci\^encias Exatas\\
Manaus, AM 69080900, Brazil}
\email{neilha@ufam.edu.br}


\begin{abstract}
In this article, we investigate the index of the MOTS given by a spatial cross section of the cosmological horizon in the
Kerr-Newman-de Sitter spacetime. We show that its index is at least one in the symmetrized sense for a small positive parameter $a$, such parameter defines the angular momentum.  Assuming a lower bound for the mass, we prove that this MOTS has index one. Also, considering an upper bound for the mass, we show that its index is at least two in the symmetrized sense. Moreover, we establish an estimate relating the area and the charge of a MOTS with index one in a Cauchy data satisfying the dominant energy condition, which give us a connection between MOTS with index one and General Relativity.
\end{abstract}

\maketitle
\section{Introduction}

In the Kerr-Newman-de Sitter spacetime $\mathcal{N}^4$, the cosmological horizon is foliated by surfaces (spatial cross sections) with vanishing outward null expansion: marginally outer trapped surface (MOTS). Furthermore, the event horizon of the Kerr-Newman-de Sitter black hole is foliated by MOTS. More generally, given a Cauchy surfaces in any spacetime, the boundary of the trapped region is an apparent horizon which is also a MOTS. Among other reasons, MOTSs are interesting objects of study because they may, in principle, locate transition zone from strong to weak gravitational field and they are frequently taken as a good sign of the location of the black hole boundary. Then, the notion of MOTS gives an effective way to understand some properties of black holes and it becomes essential to investigate in detail the geometry of MOTSs; for more details, see, e.g. \cite{BKO, CGS,HE, HL, M}.  

MOTS can be thought of as generalizations of minimal surfaces. Minimal surfaces are by definition critical points of the area functional
with respect to all possible variations. Otherwise, the stationary points of the area functional for a certain class of variations are MOTSs, more precisely, the MOTS $\Sigma$ is a smooth, orientable, closed surface in the spacetime $\mathcal{N}^4$ and the area functional is critical only for variations along the null outward pointing vector fields of $\Sigma$, see \cite{M}.
 The study of MOTS turns out to have many similarities with the study of minimal surfaces. 

The {\it{index}} of $\Sigma$ is the number of (counted with multiplicity) negative eigenvalues of the operator $L$. Also, The {\it{index}} of $\Sigma$ in the symmetrized sense is the number of (counted with multiplicity) negative eigenvalues of the symmetrized operator $L_s$. Further, we say that $\Sigma$ is {\it{stable}} if  $\lambda_1(L)\geq0$ and $\Sigma$ is {\it{unstable}} if $\lambda_1(L)<0$, where $L$ is the stability operator for MOTS (see Section \ref{Back} for definitions). It is well known the connection between stable (index zero) minimal surface and General Relativity, for example, the spatial cross sections of the event horizon in a static black holes obeying the null energy condition are stable minimal surfaces. In \cite{CLS}, the authors presented a  relation between minimal surfaces with index one and General Relativity. For instance, they proved that for a certain class of electrostatic systems, each of its unstable horizons is the solution of a one-parameter min-max problem for the area functional, in particular it has index one. In addition, they showed that the MOTS given by a  spatial cross section of the cosmological horizon in the Reissner-Nordström-de Sitter has index one.  In such model the parameter $a$ which defines the angular momentum is equal to zero whereas the parameter $a$ is positive in the Kerr-Newman-de Sitter spacetime. Motivated by these facts, we will analyze the index of the MOTS given by a spatial cross section of the cosmological horizon in the Kerr-Newman-de Sitter spacetime.

We highlight that the standard models and others physical elements fundamental for the exposition of the main results will be presented in the Section \ref{Back}. 

Taking into account the above discussion, we obtain the following result.


\begin{theo}\label{index least one}
Let $\Sigma$ be a MOTS given by a spatial cross section of the cosmological horizon in the Reissner-Nordström-de Sitter spacetime. Then
\begin{equation*}
    \lambda_1(L_s(a))<0
\end{equation*}
for a small $a>0$, and the MOTS given by a spatial cross section of the cosmological horizon in the Kerr-Newman-de Sitter spacetime has index at least one in the symmetrized sense. 
\end{theo}
As a consequence of Theorem \ref{index least one}, we obtain the following corollary (see Section \ref{Back} for definitions).
\begin{coro}\label{cor-1}
Let $\Sigma$ be a MOTS given by a spatial cross section of the cosmological horizon in the Reissner-Nordström-de Sitter spacetime. Then
\begin{equation*}
    \lambda_1(L(a))<0
\end{equation*}
for a small $a>0$, and the MOTS given by a spatial cross section of the cosmological horizon in the Kerr-Newman-de Sitter spacetime is unstable. 
\end{coro}
In light of the Theorem \ref{index least one}, it is natural to ask whether it is possible to show that the MOTS has index one. This question leads to the following result.
\begin{theo}\label{Theorem index one}
Let $\Sigma$ be a MOTS given by a spatial cross section of the cosmological horizon in the Reissner-Nordström-de Sitter spacetime
with
\begin{equation}\label{mass hypothesis}
         \frac{2\mathcal{Q}(\Sigma)^2}{3\sqrt{\beta_+}}< m,
     \end{equation}
   where $m$ is the parameter which defines the mass and  $\beta_+=\frac{3+\sqrt{9-4\Lambda\mathcal{Q}(\Sigma)^2}}{2\Lambda}$ with $\Lambda$ and $\mathcal{Q}(\Sigma)$ stand for the cosmological constant and the charge (defined in (\ref{charge expression})), respectively.
Then
\begin{equation*}
    0<\lambda_2(L_s(a)) 
\end{equation*}
for a small $a>0$, and the MOTS given by a spatial cross section of the cosmological horizon in the Kerr-Newman-de Sitter spacetime has index one in the symmetrized sense. Moreover, if the equality holds in (\ref{mass hypothesis}), then  $\Sigma$ has index one in the symmetrized sense and it is degenerate.
\end{theo}
It is also interesting to ask whether it is possible to show that the MOTS has index at least two in the symmetrized sense. For this question, we obtain the following result.

\begin{theo}\label{index least two}
Let $\Sigma$ be a MOTS given by a spatial cross section of the cosmological horizon in the Reissner-Nordström-de Sitter spacetime
with
\begin{equation}\label{mass hypothesis-2}
        m<\frac{2\mathcal{Q}(\Sigma)^2}{3\sqrt{\beta_+}},
     \end{equation}
   where $m$ is the parameter which defines the mass and  $\beta_+=\frac{3+\sqrt{9-4\Lambda\mathcal{Q}(\Sigma)^2}}{2\Lambda}$ with $\Lambda$ and $\mathcal{Q}(\Sigma)$ stand for the cosmological constant and the charge (defined in (\ref{charge expression})), respectively.
Then
\begin{equation*}
    \lambda_2(L_s(a)) <0
\end{equation*}
for a small $a>0$, and the MOTS given by a spatial cross section of the cosmological horizon in the Kerr-Newman-de Sitter spacetime has index at least two in the symmetrized sense.
\end{theo}
Now, assuming the previous setting considering the charge $\mathcal{Q}(\Sigma)$ is equal to zero, we have the following result from Theorems \ref{index least one} and \ref{Theorem index one}.

\begin{coro}\label{cor-2}
  Let $\Sigma$ be a MOTS given by a spatial cross section of the cosmological horizon in the Schwarzschild-de Sitter spacetime.
Then
\begin{equation*}
    \lambda_1(L_s(a))<0<\lambda_2(L_s(a)) 
\end{equation*}
for a small $a>0$, and the MOTS given by a spatial cross section of the cosmological horizon in the Kerr-de Sitter spacetime has index one in the symmetrized sense. Moreover,
\begin{equation*}
    r'_0<\frac{\sqrt{3}}{\sqrt{\Lambda}},
\end{equation*}
where $r'_0$ is a positive root of $\Delta_r$ in the Kerr-de Sitter case.
\end{coro}

Geometric inequalities linking the area and the physical quantities of black holes (e.g. mass, charge and angular momentum) have been attracted a large interest since the beginning of the theory; see \cite{BK,CRS,DJR,DGC, S}. Another connection between minimal surfaces with index one and General Relativity was established in \cite{CLS}, more specifically, they proved an estimate relating the area and the charge of a minimal surface with index one.

Now, we assume $(M^3,\,g,\,E)$ in a Lorentzian $4$-manifold $\mathcal{N}$ to be as part of a Cauchy data for the Einstein-Maxwell equations satisfying the dominant energy condition for the non-electromagnetic field. Then
\begin{eqnarray*}
\mu-|J|_g\geq 0,
\end{eqnarray*}
where $\mu$ is the energy density, $J$ is the momentum density of the non-electromagnetic
fields, and $(M^3,\,g)$ is an oriented Riemannian manifold.
We also take $\Sigma\subset M$ an orientable closed surface with unit normal $\nu$. The charge is expressed by 
\begin{equation}\label{charge expression}
    \mathcal{Q}(\Sigma)=\frac{1}{4\pi}\int_{\Sigma}\langle E, \nu\rangle\,d\mu,
\end{equation}
where $E\in\mathfrak{X}(M)$ is the electric field.

Inspired by the area-charge inequality showed in \cite{CLS} and the fact that the MOTS given by a spatial cross section of the cosmological horizon in the Kerr-Newman-de Sitter spacetime has index one in the symmetrized sense, we consider the previous setting to derive the next estimate in the context of MOTS instead of minimal surface which give us a connection between MOTS with index one and General Relativity (see Section \ref{Back} for definitions).
\begin{theo}\label{area-charge-result}
    Let $(M^3,\,g,\,K)$ be a maximal initial data such that
     \begin{eqnarray*}
\s_g\geq 2\Lambda+2|E|^2
\end{eqnarray*}
with positive cosmological constant $\Lambda$. Also, we consider a MOTS $\Sigma$ in $M$ with index one in the symmetrized sense such that $\Sigma$ is a topological sphere. Then
\begin{eqnarray}\label{cosmological-area-charge-ineq}
    \Lambda |\Sigma| +\frac{16\pi^2\mathcal{Q}(\Sigma)^2}{|\Sigma|}\leq 12\pi.
\end{eqnarray}
Moreover, the equality holds if and only if $ \chi_+\equiv0$, $E_{\Sigma}=c\nu$ for some constant $c$ and $(\s_g)_{\Sigma}\equiv 2\Lambda+2c^2.$
\end{theo}
\begin{remark}
    The estimates
\begin{equation*}
    \mathcal{Q}(\Sigma)^2\leq\frac{9}{4\Lambda}
\end{equation*}
and
\begin{equation*}
    4\pi\beta_{-}\leq |\Sigma|\leq 4\pi\beta_{+},
\end{equation*}
where $\beta_{\pm}=\frac{3\pm\sqrt{9-4\Lambda\mathcal{Q}^2}}{2\Lambda}$ are a straightforward application from inequality (\ref{cosmological-area-charge-ineq}).
\end{remark}

This article is organized as follows. In Section \ref{Back}, we present some definitions, the standard models used in this work, and the lemma that will be used to prove some of the main results. Finally, in Section \ref{ME}, we collect the proofs of Theorems, \ref{index least one}, \ref{Theorem index one}, \ref{index least two} and \ref{area-charge-result} and Corollaries \ref{cor-1} and \ref{cor-2}. 

\section{Background}\label{Back}

In this section, we recall some basic facts that will be essential in the proof of the main theorems.

Initially,
let $(\mathcal{N},\Bar{g})$ be a Lorentzian $4$-manifold. Consider a spacelike hypersurface $M \subset \mathcal{N}$ with unit timelike normal field $\xi$, induced metric $g$ and second fundamental form $K$, so that $(M,g,K)$ is an initial data. In our convention $K(X,Y) = \Bar{g}(\overline{\nabla}_{X}\xi,Y)$, where $\overline{\nabla}$ is the Levi-Civita connection of $(\mathcal{N},\Bar{g})$.
 
 Let $\Sigma$ be a closed surface embedded in $M$ which is two-sided, i.e, $\si$ admits a globally defined unit normal $\nu$, unique up to sign. By convention, we refer to such a choice of $\nu$ as outward pointing. Also, we denote by $h$ the metric in $\si$ induced by $g$. We define the null outward and inward pointing vector fields of $\si$ by
 $$\ell^{\pm}=\xi\pm\nu.$$

 In what follows, consider the null second fundamental form of $\Sigma\subset M$ with respect to $\ell^{\pm}$ given by
 \begin{eqnarray*}
     \chi_{\pm}(X,Y):= \Bar{g}(\overline{\nabla}_X \ell^{\pm},Y) =(K \pm A)(X,Y),
 \end{eqnarray*}
 where $A$ is the second fundamental form of $\Sigma$ inside $(M,g)$. From this, the null expansions of $\Sigma$ are expressed by
 \begin{equation*}
     \theta_{\pm}=\tr_h\chi_{\pm} =\tr_h K\pm H,
 \end{equation*}
 where $H$ is the mean curvature of $\Sigma$ in $(M,g)$. The surface $\Sigma$ is trapped if both $\theta_{+}<0$ and $\theta_{-}<0$, outer trapped if $\theta_{+}<0$, and marginally outer trapped (MOTS) if $\theta_{+}=0$. A MOTS $\Sigma$ is just a minimal surface in $M$ when $M$ is a time-symmetric, that is, $K=0$.

 Let $\Sigma_t$ denotes the one parameter family of surfaces created by deforming a MOTS $\Sigma$ an amount $t\psi$ in the outward normal direction, where $\psi:\Sigma\to\mathbb{R}$ is a smooth function. More precisely,
 $$\Sigma_t=\{\exp_x(t\psi(x)\nu_x);\, x\in\Sigma\},$$
 where $\exp$ is the exponential map of $(M,g)$. The stability of $\Sigma$ is studied by investigating
the change of $\theta_+$ along these deformations of $\Sigma$. A computation as in \cite{AMS} shows that
 \begin{eqnarray*}
     \frac{\partial\theta_+}{\partial t}\bigg\vert_{t=0}=L\psi.
 \end{eqnarray*}
 Here, the differential operator
 \begin{equation}\label{the differential operator}
     L\psi=-\Delta_h\psi+2h(Z,\nabla_h\psi)+\left(\frac{\s_h}{2}-\Ein(\ell^{+},\xi)-\frac{|\chi_{+}|^2}{2}+\diver_{h} Z-\vert Z \vert^2\right)\psi
 \end{equation}
 is the MOTS stability operator, where $\s_h$ is the scalar curvature with respect to metric $h$, $Z$ is the vector field on $\Sigma$ metrically dual to the one-form $K(\nu,\cdot)$,
 and $\Ein$ is the Einstein tensor described by
 \begin{equation*}
     \Ein=\Ric_{\Bar{g}} -\frac{\s_{\Bar{g}}}{2}\Bar{g}.
 \end{equation*}
The operator $L$ is not self-adjoint in general and its principal eigenvalue $\lambda_1(L)$ (defined as the eigenvalue with smallest real part) is real. Moreover, the associated eigenfunction $\psi_1$, $L\psi_1=\lambda_1\psi_1$, is unique up to a scale, and it can be chosen to be strictly positive; see \cite{AMS} and the references therein.
We remember that $\Sigma$ is {\it{stable}} if $\lambda_1(L)\geq0$ and $\Sigma$ is {\it{unstable}} if $\lambda_1(L)<0$.

Now, we consider the symmetrized operator $L_s: C^{\infty}(\Sigma)\longrightarrow C^{\infty}(\Sigma)$ given by
\begin{equation*}
    L_s\psi=-\Delta_h\psi+b\psi,
\end{equation*}
where $b=\Big(\frac{\s_h}{2}-\Ein(\ell^+,\xi)-\frac{|\chi_{+}|^2}{2}\Big)$. This operator is obtained by placing $Z=0$ in (\ref{the differential operator}).  We define the principal eigenvalue of the operator $L_s$ by
\begin{eqnarray*}\label{principal eigenvalue}
    \lambda_1 (L_s)=\inf_{\psi \in C^{\infty}(\Sigma)\setminus\{0\}}\frac{\int_{\Sigma}(|\nabla\psi|^2+b\psi^2)d\mu}{\int_{\Sigma}\psi^2d\mu}.
\end{eqnarray*}
The operator $L_s$ is self-adjoint, its eigenvalues can be written in increasing order, and $\lambda_1(L_s)\geq\lambda_1(L)$; see \cite{AMS,Galloway,GS,M} and the references therein.

Next, we exhibit the notion of stability in the symmetrized sense.
\begin{definition}
A MOTS $\Sigma$ is called
\begin{itemize}
    \item Stable in the symmetrized sense if $\lambda_1(L_s)\geq0$.
    \item Unstable in the symmetrized sense if $\lambda_1(L_s)<0$.
\end{itemize}
    \end{definition}
\noindent  Recall that the {\it{index}} of $\Sigma$ in the symmetrized sense is the number of (counted with multiplicity) negative eigenvalues of $L_s$. Also, we say that $\Sigma$ is degenerate if one of the eigenvalues of $L_s$ is equal to zero.

Now, we present the models and their MOTS given by a spatial cross section of the cosmological horizon for the sake of completeness. We also point out some relevant properties of these models that were used to prove some of the main results of this article.

 Let $\mathcal{N}=\mathbb{R}^2\times\mathbb{S}^2$ be the Lorentzian 4-manifold endowed with the Kerr-Newman-de Sitter metric \cite{BK, CCK} which in the Boyer-Lindquist coordinates is described by
\begin{equation}\label{Kerr-Newman-de Sitter metric}
 \Bar{g}=-\frac{\Delta_r}{\rho^2}\Big(dt-\frac{a\sin^2\theta}{\Xi}d\phi\Big)^2+\frac{\rho^2}{\Delta_r}dr^2+\frac{\rho^2}{\Delta_{\theta}}d\theta^2+\frac{\Delta_{\theta}\sin^2\theta}{\rho^2}\Big(adt-\frac{r^2+a^2}{\Xi}d\phi\Big)^2,
\end{equation}
where
\begin{align*}
\Delta_{\theta}&=1+\frac{\Lambda a^2}{3}\cos^2\theta,\\
\Delta_ r&=(r^2+a^2)\left(1-\frac{\Lambda r^2}{3}\right)-2mr+q^2,\\
\rho^2&=r^2+a^2\cos^2\theta,\\
\Xi &= 1+\frac{\Lambda a^2}{3}.
\end{align*}

\noindent The parameters $a$, $m$, and $q$ determine the angular momentum, the mass, and the charge, 
\begin{eqnarray*}
    \mathcal{J}=\frac{am}{\Xi^2},\,\,\,\mathcal{M}=\frac{m}{\Xi^2},\,\,\,\mathcal{Q}=\frac{q}{\Xi},
\end{eqnarray*}
respectively. In the sequel, we choose the values of the parameters $m>0$ and $a\geq 0$ in such a way that $\Delta_r$ admits four distinct real roots, which are noted
$$r_{--}<0<r_{-}<r_{+}<r_{c}.$$  
The roots $r_{-}$, $r_{+}$ and $r_c$ have a physical meaning, that is, $\{r=r_{-}\}$ is the inner (Cauchy) black hole horizon,  $\{r=r_{+}\}$ is the outher (Killing) black hole horizon, $\{r=r_c\}$ is the cosmological horizon, and the root $r_{--}$ has no physical meaning. The cosmological horizon appears when the cosmological constant $\Lambda$ is positive; for more details, see, e.g. \cite{BKO, BK, CLS}.

Currently, we assume $\mathcal{N}=\mathbb{R}^2\times\mathbb{S}^2$ equipped with the Reissner-Nordstr\"om-de Sitter metric given by
\begin{eqnarray*}
        \mathfrak{g}=-\frac{f(r)}{r^2}dt^2+\frac{r^2}{f(r)}dr^2+r^2g_{\mathbb{S}^2},
    \end{eqnarray*}
where $f(r)=-\frac{\Lambda}{3}r^4+r^2-2mr+\mathcal{Q}^2.$ Note that this metric is obtained by setting $a=0$ in (\ref{Kerr-Newman-de Sitter metric}). We will denote $(M,\,g,\, K)$ to be an induced initial data with $M:=\{t=0\}$ and $K$ is the second fundamental form. In the $(r,\,\theta,\,\phi)$ coordinate system, the induced metric $g$ is expressed by 
\begin{eqnarray*}
    g=\frac{r^2}{f(r)}dr^2+r^2g_{\mathbb{S}^2}.
\end{eqnarray*}
Using the change of variables
\begin{equation*}
    y=\int_0^r\frac{1}{\sqrt{f(s)}}ds,\quad \forall y.
\end{equation*}
Thus, the metric $g$ can be expressed in coordinates $(y,\, \theta,\,\phi)$ as
\begin{equation*}
g=r^2dy^2+r^2g_{\mathbb{S}^2}.
\end{equation*}
The induced metric on MOTSs given by a spatial cross sections of the horizons
\begin{equation*}
\Sigma:=\Big\{y_0;\, y_0=\int_0^{r_0}\frac{1}{\sqrt{f(s)}}ds \ \ \mbox{ with $r_0$ is a root of}  \ f(r) \Big\}   
\end{equation*}
is 
\begin{equation}\label{conformal metric}
h=r_0^2g_{\mathbb{S}^2},
\end{equation}
 where $g_{\mathbb{S}^2}$ is the metric of the round sphere with radius 1. The MOTS given by a spatial cross section of the cosmological horizon is denoted by 
 \begin{equation}\label{cosmological MOTS}
 \Sigma:=\Big\{y_c;\, y_c=\int_0^{r_c}\frac{1}{\sqrt{f(s)}}ds \ \ \mbox{ with $r_c$ is a root of}  \ f(r) \Big\}.
   \end{equation}
Considering the charge $\mathcal{Q}(\Sigma)=0$ in the Kerr-Newman-de Sitter and in the Reissner-Nordström-de Sitter models, we obtain the Kerr-de Sitter and the Schwarzschild-de Sitter models, respectively. In \cite{BLP, LP}, the authors present a nice approach relating the geometry of the spacetime models and  the General Relativity. The Reissner-Nordström-de Sitter and the Schwarzschild-de Sitter models have the parameter $a=0$ (the angular momentum is equal to zero) which implies that their MOTSs are round $2$-spheres. Otherwise, the Kerr-Newman-de Sitter and the Kerr-de Sitter models have the parameter $a>0$ (the angular momentum is positive) their MOTSs are axisymmetric oblate $2$-spheres. The MOTS introduced in (\ref{cosmological MOTS}) will be in the Kerr-Newman-de Sitter spacetime if $\mathcal{Q}(\Sigma)\neq 0$ and $a>0$, in the Reissner-Nordström-de Sitter spacetime if $\mathcal{Q}(\Sigma)\neq 0$ and $a=0$, in the Kerr-de Sitter spacetime if $\mathcal{Q}(\Sigma)= 0$ and $a>0$, and in the Schwarzschild-de Sitter spacetime if $\mathcal{Q}(\Sigma)= 0$ and $a=0$. The MOTS given by a spatial cross sections of the Cauchy and Killing horizons are stables in these models \cite{BKO, CLS}.

In the Reissner-Nordström-de Sitter spacetime $(\mathcal{N},\,\mathfrak{g})$, the cosmological constant, the mass, and the charge satisfy
\begin{equation}\label{Mass-inequality}
   0< m^2<\frac{1}{18\Lambda}\Big[1+12\mathcal{Q}^2\Lambda+(1-4\mathcal{Q}^2\Lambda)^{\frac{3}{2}}\Big].
\end{equation}
Also, $\mathcal{Q}^2\Lambda\leq\frac{1}{4}$; see \cite{B,CLS} for more details. We denote $r_{--}$, $r_{-}$, $r_+$, $r_c$; $r_1$, $r_2$, $r_3$ and $\hat{r}_1$, $\hat{r}_2$ to be the roots of 
\begin{eqnarray*}
    f(r):=-\frac{\Lambda}{3}r^4+r^2-2mr+\mathcal{Q}^2
\end{eqnarray*}
$f'(r)$ and $f''(r)$, respectively. 
\begin{lemma}\label{lemma}
The roots of $f(r),$ $f'(r)$ and $f''(r)$ have the following behavior
$$r_{--}<r_1<\hat{r}_1<0<r_{-}<r_2<r_{+}<r_3<r_{c}$$
with $\hat{r}_2\in (r_2,r_3)$ and
\begin{eqnarray*}
\left\{\begin{array}{rc}
f'(r) &\mbox{is strictly decreasing in}\,\,\, \Big(-\infty,-\sqrt{\frac{1}{2\Lambda}}\Big)\cup\Big(\sqrt{\frac{1}{2\Lambda}},+\infty\Big)\\
f'(r)  &\mbox{ is strictly increasing in}\,\,\, \Big(-\sqrt{\frac{1}{2\Lambda}},\sqrt{\frac{1}{2\Lambda}}\Big).
\end{array}\right.
\end{eqnarray*}
\end{lemma}

\begin{proof}
   In fact, we already know that the roots of $f(r)$ have the following behavior
$$r_{--}<0<r_{-}<r_{+}<r_{c}.$$
Consequently,
 $$r_{--}<r_1<0<r_{-}<r_2<r_{+}<r_3<r_c,$$
 by using Vieta's formulas and Rolle's theorem.
 Proceeding, 
$$\hat{r}_1=-\sqrt{\frac{1}{2\Lambda}}\,\,\mbox{and}\,\,\hat{r}_2=\sqrt{\frac{1}{2\Lambda}}.$$
 Thus, using Rolle's theorem again
 $$r_{--}<r_1<\hat{r}_1<0<r_{-}<r_2<r_{+}<r_3<r_c$$
and $\hat{r}_2\in (r_2,r_3).$ We also note that
\begin{eqnarray*}
\left\{\begin{array}{rc}
f''(r) &\mbox{is negative in}\,\,\, \Big(-\infty,-\sqrt{\frac{1}{2\Lambda}}\Big)\cup\Big(\sqrt{\frac{1}{2\Lambda}},+\infty\Big)\\
f''(r)  &\mbox{ is positive in}\,\,\, \Big(-\sqrt{\frac{1}{2\Lambda}},\sqrt{\frac{1}{2\Lambda}}\Big)
\end{array}\right.
\end{eqnarray*} 
which implies that
\begin{eqnarray*}
\left\{\begin{array}{rc}
f'(r) &\mbox{is strictly decreasing in}\,\,\, \Big(-\infty,-\sqrt{\frac{1}{2\Lambda}}\Big)\cup\Big(\sqrt{\frac{1}{2\Lambda}},+\infty\Big)\\
f'(r)  &\mbox{ is strictly increasing in}\,\,\, \Big(-\sqrt{\frac{1}{2\Lambda}},\sqrt{\frac{1}{2\Lambda}}\Big).
\end{array}\right.
\end{eqnarray*} 
\end{proof}

Proceeding, we present the setting that it will be used in Theorem \ref{area-charge-result}. The Einstein-Maxwell equations for uncharged matter are
  \begin{align*}
\Ric_{\Bar{g}}-\frac{\s_{\Bar{g}}}{2}\Bar{g}+\Lambda\Bar{g} &=8\pi(T_{F}+T),\\
d F &=0,\,\, \diver_{\Bar{g}} F=0. 
\end{align*}
Here $(\mathcal{N}^4,\,\Bar{g})$ is a time-oriented spacetime, $F$ is a $2$-form on $\mathcal{N}$,  $T_{F}$ is the energy-momentum tensor of the electromagnetic field given by
\begin{equation*}
    T_{F}=\frac{1}{4\pi}\Big( F\circ F-\frac{1}{4}|F|^2_{\Bar{g}}\Bar{g}\Big)
\end{equation*}
with $(F\circ F)_{ij}=\Bar{g}^{kl} F_{ik} F_{jl}$, and $T$ is a symmetric covariant $2$-tensor field on $\mathcal{N}$, which represents the energy-momentum tensor of the non-electromagnetic matter.

A Cauchy data for these equations (without a magnetic field) is a tuple
\begin{equation*}
    (M^3,\,g,\,K,\,E,\,\mu,\,J),
\end{equation*}
where $(M^3,\,g)$ is an oriented Riemannian $3$-manifold in $\mathcal{N}$ such that 
$$\mu:= T(\xi,\xi),\,\,\,\, J:=T(.,\xi)$$
are the energy density and the momentum density of the non-electromagnetic fields, respectively. Moreover, the data satisfies the constraint equations
\begin{align*}
\s_{g}+(\tr_gK)^2-|K|_g^2-2\Lambda &=2|E|^2+16\pi\mu,\\
\diver_{g}(K-(\tr_gK)_g) &=8\pi J,\\
\diver_{g}E &=0.
\end{align*} 

In this context, we say that it holds the dominant energy condition for the non-electromagnetic fields if \cite[Section 2]{DGC}
\begin{equation*}
    T(X,\,X)\geq 0,
\end{equation*}
for all future directed causal vectors. As a consequence of this condition, we obtain $\mu\geq |J|_g$ \cite[Section 8.3.4]{G}. Thus, assuming $(M,\,g,\,K)$ is maximal, that is, $\tr_gK=0$, we get 
\begin{equation*}
    \s_g\geq 2\Lambda+2|E|^2.
\end{equation*}

\section{Proof of the Main Results}\label{ME}
In this section, we shall present the proofs of Theorems, \ref{index least one}, \ref{Theorem index one}, \ref{index least two} and \ref{area-charge-result} and Corollaries \ref{cor-1} and \ref{cor-2}.

\begin{proof}[Proof of Theorem~\ref{index least one}]

By the perturbation theory for the eigenvalues in (\cite{K}, p. 63), we know that if $\lambda_1(L_s(0))<0$, then
$$\lambda_1(L_s(a))<0$$
for a small $a>0$. So, it suffices to show $\lambda_1(L_s(0))<0$ to prove the result. 

In equation (\ref{conformal metric}), we saw that the metrics $h$ and $g_{\mathbb{S}^2}$ are conformal. Then
\begin{equation*}
    \Delta_{h}\psi=\frac{1}{r_0^2}\Delta_{\mathbb{S}^2}\psi
\end{equation*}
and
\begin{equation*}
    \s_{h}=\frac{2}{r_0^2}.
\end{equation*}
Since the Reissner-Nordström-de Sitter spacetime is stationary and satisfies the null energy condition, we have $\chi_+=0$. In addition,
\begin{equation*}
     \Ein(\ell^{+},\xi)=\Lambda+\frac{\mathcal{Q}^2}{r_0^4}.
 \end{equation*}
 From this, the symmetrized operator $L_s: C^{\infty}(\Sigma)\longrightarrow C^{\infty}(\Sigma)$ follows 
\begin{eqnarray*}
  L_s\psi
  =-\Delta_{\mathbb{S}_{r_0}^2}\psi+\Big(\frac{1}{r_0^2}-\Lambda-\frac{\mathcal{Q}^2}{r_0^4}\Big)\psi,
\end{eqnarray*}
 where $\Delta_{\mathbb{S}^2_{r_0}}$ stands for the Laplacian operator of the round sphere with radius $r_0$. Now,

 \begin{equation*}
   \Delta_h=  \Delta_{\mathbb{S}^2_{r_0}}=r_0^{-2}\Delta_{\mathbb{S}^2_1}
 \end{equation*}
 and
 \begin{eqnarray*}
     -\Delta_{\mathbb{S}^2_{r_0}}\psi=\Bar{\lambda}\psi
     \Leftrightarrow r_0^{-2}\Delta_{\mathbb{S}_1^2}\psi
     =-\Bar{\lambda}\psi\Leftrightarrow r_0^{-2}k(k+1)\psi=\Bar{\lambda}\psi,
 \end{eqnarray*}
 where $\hat{\lambda}_{k+1}=k(k+1)$ is the spectrum of $-\Delta_{\mathbb{S}_1^2}$ and $\Bar{\lambda}_{k+1}=r_0^{-2}k(k+1)$ is the spectrum of $-\Delta_{\mathbb{S}^2_{r_0}}$. This leads us to
 \begin{eqnarray*}
     -\Delta_{\mathbb{S}^2_{r_0}}\psi+b\psi=\lambda\psi\Leftrightarrow -\Delta_{\mathbb{S}^2_{r_0}}\psi=(\lambda-b)\psi=\Bar{\lambda}\psi,
 \end{eqnarray*}
 where
 $$b=\frac{1}{r_0^2}-\Lambda-\frac{\mathcal{Q}^2}{r_0^4}.$$
 Hence, the spectrum of $L_s(0)$ is given by
 \begin{equation}\label{spectrum of L_s RNdS}
     \lambda_{k+1}(L_s(0))=r_0^{-2}k(k+1)+\frac{1}{r_0^2}-\Lambda-\frac{\mathcal{Q}^2}{r_0^4}.
 \end{equation}
 Proceeding, we get
 \begin{equation*}
\lambda_1(L_s(0))=\frac{1}{r_0^4}(-\Lambda r_0^4+r_0^2-\mathcal{Q}^2).
 \end{equation*}
 Taking $x=r_0^2$, we obtain
 $$u(x)=-\Lambda x^2+x-\mathcal{Q}^2.$$
 Thus,
 \begin{eqnarray*}
\left\{\begin{array}{rc}
u(x)<0 &\mbox{if}\,\,\,x\in (0,\alpha_{-})\cup(\alpha_{+},+\infty)\\
u(x)>0  &\mbox{ if}\,\,\,x\in (\alpha_{-},\alpha_{+}),
\end{array}\right.
\end{eqnarray*}
where $\alpha_{-}=\frac{1-\sqrt{1-4\Lambda\mathcal{Q}^2}}{2\Lambda}$ and $\alpha_{+}=\frac{1+\sqrt{1-4\Lambda\mathcal{Q}^2}}{2\Lambda}$.
Consequently,
 \begin{eqnarray*}
\left\{\begin{array}{rc}
\lambda_1(L_s(0))<0 &\mbox{if}\,\,\,r_0\in (0,\sqrt{\alpha_{-}})\cup(\sqrt{\alpha_{+}},+\infty)\\
\lambda_{1}(L_s(0))>0  &\mbox{ if}\,\,\,r_0\in (\sqrt{\alpha_{-}},\sqrt{\alpha_{+}}).
\end{array}\right.
\end{eqnarray*}
We know that
\begin{eqnarray*}
    f'(\sqrt{\alpha_{+}})=
    2\Bigg[\frac{\sqrt{\alpha_{+}}}{3}(2-\sqrt{1-4\mathcal{Q}^2\Lambda})-m\Bigg].
\end{eqnarray*}
Using (\ref{Mass-inequality}) and the fact that
\begin{equation*}
    \frac{1}{18\Lambda}\Big[1+12\mathcal{Q}^2\Lambda+(1-4\mathcal{Q}^2\Lambda)^{\frac{3}{2}}\Big]=\Bigg[\frac{\sqrt{\alpha_{+}}}{3}(2-\sqrt{1-4\mathcal{Q}^2\Lambda})\Bigg]^2,
\end{equation*}
we achieve
\begin{eqnarray*}
    f'(\sqrt{\alpha_{+}})>0.
\end{eqnarray*}
In Lemma \ref{lemma}, we saw that
$$r_3,\,\,\,\sqrt{\alpha_{+}},\,\,\, r_c\in\Big(\sqrt{\frac{1}{2\Lambda}},+\infty\Big)$$
with $r_3<r_c$, and $f'(r)$ is strictly decreasing in $\Big(\sqrt{\frac{1}{2\Lambda}},+\infty\Big)$. Then,
$$\sqrt{\alpha_{+}}<r_3<r_c.$$
Consequently,
$$\lambda_1(L_s(0))=\frac{1}{r_0^4}(-\Lambda r_0^4+r_0^2-\mathcal{Q}^2)<0$$
when $r_0=r_c$. Therefore, the eigenvalue $\lambda_1(L_s(a))$ is negative and the MOTS given by a spatial cross section of the
cosmological horizon in the Kerr-Newman-de Sitter spacetime has index at least one in
the symmetrized sense.
\end{proof}
\begin{remark}
The argument that $f(r)$ is positive in the region where $r\in(r_+,r_c)$ and negative in the region where $r\in(r_{-},r_{+})\cup(r_c,+\infty)$ can be used to give a different proof for Theorem \ref{index least one}.
 \end{remark}
\begin{proof}[Proof of Corollary~\ref{cor-1}]
In Section \ref{Back}, we saw that
\begin{equation*}
    \lambda_1(L(a))\leq\lambda_1(L_s(a))
\end{equation*}
and by Theorem \ref{index least one}, we know that
\begin{equation*}
    \lambda_1(L_s(a))<0.
\end{equation*}
Hence,
\begin{equation*}
     \lambda_1(L(a))<0
\end{equation*}
for small $a>0$. Then, the MOTS given by a spatial cross section of the cosmological horizon in the Kerr-Newman-de Sitter spacetime is unstable.
\end{proof}
\begin{proof}[Proof of Theorem~\ref{Theorem index one}]
By using the perturbation theory for the eigenvalues in (\cite{K}, p. 63), we have that if $\lambda_2(L_s(0))>0$, then
$$\lambda_2(L_s(a))>0$$
for a small $a>0$. From this and Theorem \ref{index least one}, we conclude that
$$\lambda_1(L_s(a))<0<\lambda_2(L_s(a))$$
for a small $a>0$. Then it is enough to show $\lambda_2(L_s(0))>0$ to prove that the MOTS given by a spatial cross section of the cosmological horizon in the Kerr-Newman-de Sitter spacetime has index one in the symmetrized sense.

Proceeding, we deduce that
 \begin{equation*}
     \lambda_2(L_s(0))=\frac{1}{r_0^4}(-\Lambda r_0^4+3r_0^2-\mathcal{Q}^2)
 \end{equation*}
from (\ref{spectrum of L_s RNdS}). If we assume $x=r_0^2$, we get
\begin{equation*}
    v(x)=-\Lambda x^2+3x-\mathcal{Q}^2.
\end{equation*}
So,
\begin{eqnarray*}
\left\{\begin{array}{rc}
v(x)<0 &\mbox{if}\,\,\,x\in (0,\beta_{-})\cup(\beta_{+},+\infty)\\
v(x)>0  &\mbox{ if}\,\,\,x\in (\beta_{-},\beta_{+}),
\end{array}\right.
\end{eqnarray*}
where $\beta_{-}=\frac{3-\sqrt{9-4\Lambda\mathcal{Q}^2}}{2\Lambda}$ and $\beta_{+}=\frac{3+\sqrt{9-4\Lambda\mathcal{Q}^2}}{2\Lambda}$. Thus, we conclude that
 \begin{eqnarray*}
\left\{\begin{array}{rc}
\lambda_2(L_s(0))<0 &\mbox{if}\,\,\,r_0\in (0,\sqrt{\beta_{-}})\cup(\sqrt{\beta_{+}},+\infty)\\
\lambda_{2}(L_s(0))>0  &\mbox{ if}\,\,\,r_0\in (\sqrt{\beta_{-}},\sqrt{\beta_{+}}).
\end{array}\right.
\end{eqnarray*}
Furthermore,
 \begin{eqnarray*}
f(\sqrt{\alpha_{+}})
&=&\frac{2\sqrt{\alpha_{+}}[(1+\sqrt{1-4\Lambda\mathcal{Q}^2})(5-\sqrt{1-4\Lambda\mathcal{Q}^2})+12\Lambda\mathcal{Q}^2]}{24\Lambda\sqrt{\alpha_{+}}}
-2m\sqrt{\alpha_{+}}.
\end{eqnarray*}
Using that
\begin{eqnarray}\label{first equality}
    \Bigg[\frac{1+12\Lambda\mathcal{Q}^2+(1-4\Lambda\mathcal{Q}^2)^{\frac{3}{2}}}{18\Lambda}\Bigg]^{\frac{1}{2}}
  =  \frac{(1+\sqrt{1-4\Lambda\mathcal{Q}^2})(5-\sqrt{1-4\Lambda\mathcal{Q}^2})+12\Lambda\mathcal{Q}^2}{24\Lambda \sqrt{\alpha_{+}}}
 \end{eqnarray}
and (\ref{Mass-inequality}), we have
$$f(\sqrt{\alpha_{+}})>0.$$
Also,
 \begin{equation*}
       \sqrt{\alpha_{+}}\in(\sqrt{\beta_{-}},\,\sqrt{\beta_{+}}).
    \end{equation*}
    In another direction,
   \begin{eqnarray*}
f(\sqrt{\beta_{+}})&=&(\sqrt{\beta_{+}})^2[1-\frac{\Lambda}{3}(\sqrt{\beta_{+}})^2]-2m\sqrt{\beta_{+}}+\mathcal{Q}^2\\
&=&\Big(\Big(\frac{3+\sqrt{9-4\Lambda\mathcal{Q}^2}}{2\Lambda}\Big)^{\frac{1}{2}}\Big)^2\Big[1-\frac{\Lambda}{3}\Big(\Big(\frac{3+\sqrt{9-4\Lambda\mathcal{Q}^2}}{2\Lambda}\Big)^{\frac{1}{2}}\Big)^2\Big]\\
&-&2m\Big(\frac{3+\sqrt{9-4\Lambda\mathcal{Q}^2}}{2\Lambda}\Big)^{\frac{1}{2}}+\mathcal{Q}^2\\
&=&\frac{4\mathcal{Q}^2}{3}-2m\sqrt{\beta_{+}}.
\end{eqnarray*}
This combined with (\ref{mass hypothesis}) gives
    \begin{equation*}
        f(\sqrt{\beta_{+}})<0.
    \end{equation*}
     From the above discussion, there exist a root of $f(r)$ in
     \begin{equation*}
         I:=(\sqrt{\alpha_{+}},\,\sqrt{\beta_{+}}).
     \end{equation*}
     In Lemma \ref{lemma}, we saw that
     $$r_2<r_+<r_3<r_c,$$
     $$\sqrt{\alpha_{+}}<r_3<r_c,$$
     and $\hat{r}_2=\sqrt{\frac{1}{2\Lambda}}\in (r_2,r_3)$. This implies $r_+\in I$ or $r_c\in I$ or $r_+,r_c\in I$. If only $r_+\in I$, we have at least two roots bigger than $r_+$ for $f(r)$. In fact, taking into account that
$$f(\sqrt{\beta_{+}})<0,$$ 
so there is $\sqrt{\beta_{+}}<\Bar{r}$ such that $f(\Bar{r})>0$, otherwise $r_c$ will be a local maximum point and $f'(r_c)=0$ which leads to a contradiction with the fact that the roots of $f(r)$ and $f'(r)$ are distinct. Further,
$$\lim\limits_{r \to {+\infty}}f(r)=-\infty.$$
Thereby, there is one root in
$$(\sqrt{\beta_{+}},\,\Bar{r})$$
and another root in $(\Bar{r},+\infty)$ which leads to a contradiction with the fact that $f(r)$ has four roots and its biggest root is $r_c$. Thus, $r_c\in I$ and
$$\lambda_2(L_s(0))=\frac{1}{r_0^4}(-\Lambda r_0^4+3r_0^2-\mathcal{Q}^2)>0$$
when $r_0=r_c$. 

Now, assuming 
     \begin{equation*}
       m=  \frac{2\mathcal{Q}^2}{3\sqrt{\beta_{+}}},
     \end{equation*}
     we have $\sqrt{\beta_{+}}$ is a root of $f(r)$. Furthermore,
     
      \begin{eqnarray*}
         r_+<\sqrt{\alpha_{+}}<\sqrt{\beta_{+}},
     \end{eqnarray*}
     since $f(r)$ is positive in the region where $r\in(r_+,r_c)$ and it is negative in the region where $r\in(r_{-},r_{+})\cup(r_c,+\infty)$ what shows that
     \begin{equation*}
         r_+\neq\sqrt{\beta_{+}}.
     \end{equation*}
Then, $r_c=\sqrt{\beta_{+}}$, $\lambda_1(L_s(0))<0$, and $\lambda_2(L_s(0))=0$ which implies that $\Sigma$ has index one and it is degenerate.
\end{proof}

\begin{remark}
Note that
 \begin{eqnarray*}
     \frac{2\mathcal{Q}^2}{3\sqrt{\beta_{+}}}<\Bigg[\frac{1+12\Lambda\mathcal{Q}^2+(1-4\Lambda\mathcal{Q}^2)^{\frac{3}{2}}}{18\Lambda}\Bigg]^{\frac{1}{2}}.
 \end{eqnarray*}
In fact, by using (\ref{first equality}), the fact that
 \begin{eqnarray*}
\frac{(1+\sqrt{1-4\Lambda\mathcal{Q}^2})(5-\sqrt{1-4\Lambda\mathcal{Q}^2})+12\Lambda\mathcal{Q}^2}{24\Lambda}&=&
\frac{1+\sqrt{1-4\Lambda\mathcal{Q}^2}+4\Lambda\mathcal{Q}^2}{6\Lambda}\\
&>&\frac{2\mathcal{Q}^2}{3}
 \end{eqnarray*}
 and
 \begin{eqnarray*}
    \frac{1}{\sqrt{\alpha_{+}}}>\frac{1}{\sqrt{\beta_{+}}}
 \end{eqnarray*}
 the asserted inequality holds. This shows that it is reasonable to ask the hypothesis (\ref{mass hypothesis})
 in the Theorem \ref{Theorem index one} since we know that the parameter $m$ satisfies (\ref{Mass-inequality}).
 \end{remark}

 \begin{proof}[Proof of Theorem~\ref{index least two}]
From the perturbation theory for the eigenvalues in (\cite{K}, p. 63), we already know that if $\lambda_2(L_s(0))<0$, then
$$\lambda_2(L_s(a))<0$$
for a small $a>0$. Thus, it is enough to prove $\lambda_2(L_s(0))<0$ to show the Theorem \ref{index least two}. In Theorem \ref{Theorem index one}, we saw 

 \begin{eqnarray*}
\left\{\begin{array}{rc}
\lambda_2(L_s(0))<0 &\mbox{if}\,\,\,r_0\in (0,\sqrt{\beta_{-}})\cup(\sqrt{\beta_{+}},+\infty)\\
\lambda_{2}(L_s(0))>0  &\mbox{ if}\,\,\,r_0\in (\sqrt{\beta_{-}},\sqrt{\beta_{+}}),
\end{array}\right.
\end{eqnarray*}
where $\beta_{-}=\frac{3-\sqrt{9-4\Lambda\mathcal{Q}^2}}{2\Lambda}$ and $\beta_{+}=\frac{3+\sqrt{9-4\Lambda\mathcal{Q}^2}}{2\Lambda}$. We also saw that
\begin{equation*}
    f(\beta_{+})=\frac{4\mathcal{Q}^2}{3}-2m\sqrt{\beta_{+}}.
\end{equation*}
By (\ref{mass hypothesis-2}), we obtain
     \begin{equation*}
    f(\beta_{+})>0.
\end{equation*}
Using
$$\lim\limits_{r \to {+\infty}}f(r)=-\infty,$$
so there is $\sqrt{\beta_+}<\Tilde{r}$ such that $f(\Tilde{r})<0$. Therefore, we have at least one root of $f(r)$ in $(\sqrt{\beta_{+}},\, \Tilde{r})$. Taking into account that $f(r)$ is positive in the region where $r\in (r_{+},\,r_c)$ and it is negative in the region where $r\in(r_{-},r_{+})\cup(r_c,+\infty)$, we get $\sqrt{\beta_{+}}\in (r_{+},\,r_c)$ since $f(\sqrt{\beta_{+}})>0$. Thereby, the root $r_c\in(\sqrt{\beta_{+}},\,\Tilde{r})$ and $\lambda_2(L_s(0))<0$. Then, the eigenvalue $\lambda_2(L_s(a))$ is negative and the MOTS given by a spatial cross section of the cosmological horizon in the Kerr-Newman-de Sitter spacetime has index at least two in
the symmetrized sense.
 \end{proof}

Next, we show a simpler proof to obtain the sign of the second eigenvalue of the symmetrized operator $L_s$ in the Kerr-de Sitter
spacetime.
\begin{proof}[Proof of Corollary~\ref{cor-2}]
By Theorem \ref{index least one}
\begin{equation*}
    \lambda_1(L_s(a))<0
\end{equation*}
for a small $a>0$. Also, the perturbation theory for the eigenvalues in (\cite{K}, p. 63) guarantees that if $\lambda_2(L_s(0))>0$, then
$$\lambda_2(L_s(a))>0$$
for a small $a>0$. So, it suffices to show $\lambda_2(L_s(0))>0$ to prove that the MOTS given by a spatial cross section of the cosmological horizon in the Kerr-de Sitter spacetime has index one in the symmetrized sense.

Now, we use (\ref{spectrum of L_s RNdS}) with $\mathcal{Q}=0$ to obtain
 \begin{equation*}
\lambda_2(L_s(0))=\frac{3}{r_0^2}-\Lambda.
 \end{equation*}
 Since $m>0$ and $r_0$ is a positive root of 
\begin{equation*}
    0=r^2\Big(1-\frac{\Lambda r^2}{3}\Big)-2mr_,
\end{equation*}
we have
$$\Big(1-\frac{\Lambda r_0^2}{3}\Big)>0.$$
Multiplying the above inequality by $\frac{3}{r_0^2}$, we achieve
 \begin{eqnarray*}
     \frac{3}{r_0^2}-\Lambda>0.
 \end{eqnarray*}
 So,
 $$\lambda_2(L_s(0))>0.$$
  
In the Kerr-de Sitter case, we can also see that if $r'_0$ is a positive root of 

$$0=(r^2+a^2)\left(1-\frac{\Lambda r^2}{3}\right)-2mr,$$
we get
$$1-\frac{\Lambda r'_{0}{^2}}{3}>0$$
and

\begin{equation*}
r'_0<\frac{\sqrt{3}}{\sqrt{\Lambda}}.
\end{equation*}

\end{proof}


\begin{proof}[Proof of Theorem~\ref{area-charge-result}]
    The proof consists in using the Hersch's trick. Let $\psi_1$ be the first eigenfunction of the symmetrized operator 
     given by
    \begin{equation*}
    L_s=-\Delta_h+\frac{\s_h}{2}-\Ein(\ell^+,\xi)-\frac{|\chi_+|^2}{2},
\end{equation*}
which we can choose to be positive. Let $\Phi$ be a conformal map from $\Sigma$ to $\mathbb{S}^2\subset\mathbb{R}^3$. So, we can take the integral 
\begin{eqnarray*}
    \int_{\Sigma}(\Psi\circ\Phi)\psi_1\,d\mu\in\mathbb{R}^3,
\end{eqnarray*}
where $\Psi$ is a M\"obius transformation of $\mathbb{S}^2$. Using the fact that $\psi_1$ is positive, we can choose $\Psi$ in such a way that the previous integral vanishes; see \cite{LY}.

The coordinate $(\Psi\circ\Phi)_i=\Tilde{\psi}_i$ of $(\Psi\circ\Phi)$ is orthogonal to $\psi_1$. Since $\Sigma$ has index one, we have
\begin{eqnarray*}
    \int_{\Sigma}\Big(|\nabla_h\Tilde{\psi_i}|^2+\Big(\frac{\s_h}{2}-\Ein(\ell^+,\xi)-\frac{|\chi_+|^2}{2}\Big)\Tilde{\psi}_i^2\Big)d\mu\geq 0,\,\,i=1,2,3.
\end{eqnarray*}
Summing these three inequalities, we get

\begin{eqnarray}\label{summing-inequality}
    0\leq\int_{\Sigma}\Big(|\nabla_h(\Psi\circ\Phi)|^2+\Big(\frac{\s_h}{2}-\Ein(\ell^+,\xi)-\frac{|\chi_+|^2}{2}\Big)\Big)d\mu.
\end{eqnarray}
By \cite{YY}, we have
\begin{equation*}
    \int_{\Sigma}|\nabla_h(\Psi\circ\Phi)|^2d\mu=8\pi deg (\Psi\circ\Phi)
\end{equation*}
and
$$\Ein(\ell^+,\xi)\geq \Lambda +|E|^2,$$
where we use 
$$\s_g\geq 2\Lambda+2|E|^2$$
together with Einstein equations. 

The estimate (\ref{summing-inequality}) becomes

\begin{eqnarray}\label{last inequality}
    0&\leq& 8\pi deg (\Psi\circ\Phi) +\int_{\Sigma}\Big(\frac{\s_h}{2}-\Lambda-|E|^2\Big)d\mu\nonumber\\
    &=& 8\pi deg(\Phi)+4\pi-\Lambda |\Sigma|-\int_{\Sigma}|E|^2d\mu,
\end{eqnarray}
where we use the Gauss-Bonnet Theorem and $g(\Sigma)=0$ to obtain (\ref{last inequality}).

By \cite{RR}, we can choose $\Phi$ such that
\begin{equation*}
    deg(\Phi)\leq 1+\Big[\frac{g(\Sigma)+1}{2}\Big]
\end{equation*}
and using the Cauchy-Schwarz inequality in (\ref{last inequality}), we get
\begin{eqnarray*}
    0&\leq& 12\pi-\Lambda|\Sigma|-\int_{\Sigma} \langle E, \nu\rangle^2d\mu \\
    &\leq& 12\pi -\Lambda |\Sigma|-\frac{1}{|\Sigma|}\Big(\int_{\Sigma} \langle E, \nu\rangle d\mu\Big)^2,
\end{eqnarray*}
where we used the H\"older's inequality in the previous estimate. Also, taking into account (\ref{charge expression}),
we achieve
\begin{equation*}
    0\leq 12\pi-\Lambda |\Sigma| -\frac{16\pi^2\mathcal{Q}(\Sigma)^2}{|\Sigma|}.
\end{equation*}
If equality holds in the above estimate, then we have equality in all previous inequalities. Hence,
\begin{equation*}
    \chi_+\equiv 0.
\end{equation*}
At the same time, the equality on the Cauchy-Schwarz inequality implies $|E|^2=\langle E, \nu\rangle^2$ and $E=c\nu$ for some function $c=\langle E, \nu\rangle$. In another direction, the equality on H\"older's estimate gives us that $\langle E, \nu\rangle$ is constant. So,

$$(\s_g)_{\Sigma}\equiv 2\Lambda+2c^2.$$

\end{proof}

 \noindent
{\bf Acknowledgments}. I would like to thank Detang Zhou, Lucas Ambrozio, Vanderson Lima, Ernani Ribeiro Jr., and  Alcides de Carvalho Jr. for their valuable suggestions in this work. I would like to express my gratitude to Professors Patricia Guidolin, Evandro Manica, Roberto Crist\'ov\~ao, and Dagoberto Rizzotto for their support. I also want to thank the Instituto de Matem\'atica e Estatística of Universidade Federal do Rio Grande do Sul (UFRGS) for the hospitality, where the research in this paper was developed, and I was partially supported by Funda\c c\~ao Arthur Bernardes-Instituto Serrapilheira Professor Vanderson Lima grant ``How to detect black holes via geometric inequalities?''.

\end{document}